\documentclass[12pt,a4paper]{amsart}

\usepackage{amssymb,latexsym,graphicx}
\usepackage{color}
\usepackage{colortbl}
\usepackage{subcaption, float}
\usepackage[cp1252]{inputenc} 
\usepackage[T1]{fontenc}
\usepackage{lmodern}
\usepackage[english]{babel}
\usepackage{enumerate}
\usepackage[section]{placeins}
\usepackage[pagebackref=true]{hyperref}
\hypersetup{pdftex,colorlinks=true,linkcolor=blue,citecolor=blue,urlcolor=blue}

\usepackage{currfile}

\usepackage[margin=3cm]{geometry}

\newtheorem{theorem}{Theorem}[section]
\newtheorem{lemma}[theorem]{Lemma}

\newtheorem{corollary}[theorem]{Corollary}

\newtheorem{remark}[theorem]{Remark}
\newtheorem{remarks}[theorem]{Remarks}

\numberwithin{equation}{section}
\numberwithin{figure}{section}
\numberwithin{table}{section}

\definecolor{purple}{RGB}{127,0,255}

\definecolor{lgray}{gray}{0.90}

\newcommand{\C}{\mathbb{C}}

\newcommand{\K}{\mathbb{K}}
\newcommand{\N}{\mathbb{N}}
\newcommand{\Nb}{\mathbb{N^{\bullet}}}

\newcommand{\R}{\mathbb{R}}
\newcommand{\bS}{\mathbb{S}}
\newcommand{\T}{\mathbb{T}}
\newcommand{\Z}{\mathbb{Z}}

\newcommand{\cC}{\mathcal{C}}

\newcommand{\cE}{\mathcal{E}}

\newcommand{\cK}{\mathcal{K}}
\newcommand{\cL}{\mathcal{L}}

\newcommand{\cT}{\mathcal{T}}

\newcommand{\lambdah}{\hat{\lambda}}

\newcommand{\set}[1]{\left\lbrace #1 \right\rbrace}
\newcommand{\pib}{\frac{\pi}{2}}

\newcommand{\pid}{\frac{\pi}{4}}

\newcommand{\pf}{\emph{Proof.~}}

\title[Courant-sharp eigenvalues]{Courant-sharp eigenvalues of compact flat surfaces: Klein bottles and cylinders}

\date{\today}

\author[P. B\'{e}rard]{Pierre B\'erard}
\author[B. Helffer]{Bernard Helffer}
\author[R.Kiwan]{Rola Kiwan}

\address{PB: Universit\'{e} Grenoble Alpes and CNRS\\
Institut Fourier, CS 40700\\
F38058 Grenoble Cedex 9, France.}
\email{pierrehberard@gmail.com}

\address{BH: Laboratoire Jean Leray, Uni\-ver\-si\-t\'{e} de Nantes and CNRS\\
F44322 Nantes Cedex, France and LMO (Uni\-versi\-t\'e Paris-Sud).}
\email{Bernard.Helffer@univ-nantes.fr}

\address{RK: American University in Dubai, P.O.Box 28282, Dubai, United Arab Emirates.}
\email{rkiwan@aud.edu}

\thanks{The authors would like to thank the referee for his comments.}



\date{\today~(\currfilename)}

\begin{document}

\begin{abstract}{The question of determining for which eigenvalues there exists an eigenfunction which has the same number of nodal domains as the label of the associated eigenvalue (Courant-sharp property) was motivated by the analysis of minimal spectral partitions. In previous works, many examples have been analyzed corresponding to squares, rectangles, disks, triangles, tori, M\"{o}bius strips,\ldots . A natural toy model for further investigations is the flat Klein bottle, a non-orientable surface with Euler characteristic $0$, and particularly the Klein bottle associated with the square torus, whose eigenvalues have higher multiplicities.  In this note, we prove that the only Courant-sharp eigenvalues of the flat Klein bottle associated with the square torus (resp. with square fundamental domain) are the first and second eigenvalues. We also consider the flat cylinders $(0,\pi) \times \bS^1_r$ where $r \in \set{0.5,1}$ is the radius of the circle $\bS^1_r$, and we show that the only Courant-sharp Dirichlet eigenvalues of these cylinders are the first and second eigenvalues.}

\end{abstract}

\maketitle

\section{Introduction}\label{S-int}

Given a compact Riemannian surface $(M,g)$, we write the eigenvalues of the Laplace-Beltrami operator $-\Delta_g$,
\begin{equation}\label{E-int-2}
\lambda_1(M) < \lambda_2(M) \le \lambda_3(M) \le \ldots \,,
\end{equation}
in nondecreasing order, starting from the label $1$, with multiplicities accounted for. If the boundary $\partial M$ of $M$ is non-empty, we consider Dirichlet eigenvalues.\medskip

Courant's nodal domain theorem (1923) states that  any  eigenfunction associated with the eigenvalue $\lambda_k$ has at most $k$ nodal domains (connected components of the complement of the zero set of $u$).   The eigenvalue $\lambda_k$ is called \emph{Courant-sharp} if there exists an associated eigenfunction with precisely $k$ nodal domains. It follows from Courant's theorem that the eigenvalues $\lambda_1$ and $\lambda_2$ are Courant-sharp, and that $\lambda_{k-1} < \lambda_k$ whenever $\lambda_k$ is Courant-sharp. Courant-sharp eigenvalues appear naturally in the context of partitions, \cite{HHOT}. \medskip

Back in 1956, Pleijel proved that there are only finitely many Courant-sharp eigenvalues. More precisely, Pleijel's original proof \cite{Pl} applied to Dirichlet eigenvalues of bounded domains in $\R^2$. It was later adapted to more general domains (\cite{Pe}) and to general closed Riemannian manifolds (\cite{BM}). \medskip

Courant-sharp eigenvalues of flat tori and M\"{o}bius strips were studied in \cite{HHO3,Le,BH2016} and \cite{BHK2020} respectively. In view of the classification of complete, flat surfaces (see \cite[p.~222-223]{KN1-1963}, or \cite[Chap.~2.5]{Wolf-2011}) it is natural to investigate the Courant-sharp property for the other compact flat surfaces, Klein bottles and cylinders. The purpose of the present paper is to prove the following theorems.

\begin{theorem}\label{T-int-2K}
Let $\K_1$ denote the flat Klein bottle associated with the square torus, resp. $\K_2$ the flat Klein bottle whose fundamental domain is a square. Then the only Courant-sharp eigenvalues of $\K_c, c\in \set{1,2}$, are $\lambda_1$ and $\lambda_2$.
\end{theorem}%

\begin{theorem}\label{T-int-2C}
Let $\cC_r$ denote the flat cylinder $(0,\pi)\times \bS^1_r$, with the product metric. Here $\bS^1_r$ denotes the circle with radius $r$. Then, for $r \in \set{\frac 12, 1}$, the only Courant-sharp Dirichlet eigenvalues of $\cC_r$, are $\lambda_1$ and $\lambda_2$.
\end{theorem}%

\begin{figure}[!ht]
\begin{center}
\includegraphics[width=10cm]{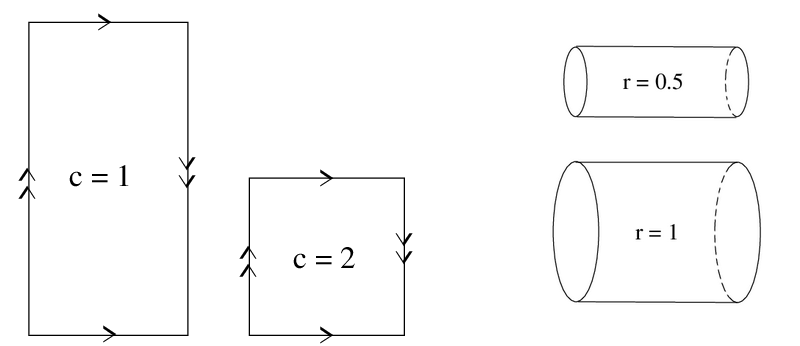}
\caption{The Klein bottles $\K_c$ (left) and the cylinders $\cC_r$ (right).} \label{F-int-K}
\end{center}
\end{figure}

Courant-sharp eigenvalues have previously been determined for several compact surfaces. We refer to the following papers and their bibliographies.
\begin{itemize}
\item[$\diamond$] Closed surfaces: round 2-sphere and projective plane, \cite{Ley96}; flat tori, \cite{HHO3, Le, BH2016}.
\item[$\diamond$] Plane domains and compact surfaces with boundary (different boundary conditions might be considered): square, \cite{Pl, BH2015, HP2015, GH1, GH2}; equilateral triangle, \cite{BH2016}; $2$-rep-tiles, \cite{BaBeFa2017}; thin cylinders, \cite{HHO2}; M\"{o}bius strips, \cite{BHK2020}.
\end{itemize}

There are also a few results in higher dimensions, see for example \cite{HPS2016, Le2016, HK2017}. \medskip

Most of the papers mentioned above adapt the method introduced by Pleijel in \cite{Pl} to the example at hand.
\medskip

The paper is organized as follows. In Section~\ref{S-ple}, we recall the main lines of Pleijel's method. More precisely, we show how a lower bound on the ratio $\frac{\lambda_k(M)}{k}$ and a lower bound on the Weyl counting function can be used to restrict the search for Courant-sharp eigenvalues to a finite set of eigenvalues. To actually determine the Courant-sharp eigenvalues one then needs to analyze the nodal patterns of eigenfunctions in a finite set of eigenspaces. In Section~\ref{S-pre}, we recall some basic facts concerning Klein bottles. In Sections~\ref{S-kle} and \ref{S-smal}, we adapt Pleijel's method (Section~\ref{S-ple}) to the flat Klein bottles  $\K_1$ and $\K_2$. In Section~\ref{S-cyl}, we adapt Pleijel's method (Section~\ref{S-ple}) to the flat cylinders $\cC_r$, with $r \in \set{\frac 12, 1}$.

\section{Pleijel's method summarized}\label{S-ple}

In order to prove that the number of Courant-sharp (Dirichlet) eigenvalues of a compact Riemannian surface $(M,g)$ is finite, one only needs two ingredients,
\begin{enumerate}
\item the classical Weyl asymptotic law,  $W_M(\lambda) = \frac{|M|}{4 \pi}\,
\lambda + O(\sqrt{\lambda})$, for the counting function $W_M(\lambda) = \set{j ~|~ \lambda_j(M) < \lambda}$, where $|M|$ denotes the area of $M$;
\item an inequality \`{a} la Faber-Krahn for domains $\omega \subset M$
with small enough area,
\begin{equation}\label{E-ple-2}
\delta_1(\omega) \ge (1-\varepsilon)^2\, \frac{\pi j_{0,1}^2}{|\omega|} \text{~~provided that~~} |\omega| \le C(M,\varepsilon)\, |M|\,,
\end{equation}
for some $\varepsilon \in [0,1)$, and for some constant $C(M,\varepsilon)$. Here, $\delta_1(\omega)$ denotes the least Dirichlet eigenvalue of $\omega$, and $j_{0,1}$ is the first positive zero of the Bessel function $J_0$ ($j_{0,1} \approx 2.404825$).
\end{enumerate}

Note that the right-hand side of \eqref{E-ple-2} is equal to $(1-\varepsilon)^2\, \delta_1(\omega^*)$, where $\omega^*$ denotes a disk in $\R^2$ with area equal to $|\omega|$. The existence of such an inequality (for a given $\varepsilon \in (0,1)$) follows from the asymptotic isoperimetric inequality proved in \cite[Appendice~C]{BM}, and from the usual symmetrization argument to translate the isoperimetric inequality into a Faber-Krahn type inequality for the first Dirichlet eigenvalue, \cite[Lemme~15]{BM}.\medskip

In order to give more quantitative information on the Courant-sharp eigenvalues, one needs a lower bound on the Weyl function, in the form,
\begin{equation}\label{E-ple-4}
W_M(\lambda) \ge \frac{|M|}{4 \pi}\, \lambda - A(M)\, \sqrt{\lambda} - B(M)\,, \hspace{3mm} \text{for all~} \lambda \ge 0
\end{equation}
for some constants $A(M)$ and $B(M)$ depending on the geometry of $(M,g)$.

\begin{lemma}\label{L-ple-2}
Assume that $(M,g)$ satisfies \eqref{E-ple-2}, with $\varepsilon < 1 - \frac{2}{j_{0,1}}$, and \eqref{E-ple-4}. Let $\lambda_k(M)$ be a Courant-sharp eigenvalue with $k \ge \frac{1}{C(M,\varepsilon)}$. Then,
\begin{align}
\frac{\lambda_k(M)}{k} &\ge (1-\varepsilon)^2 \frac{\pi j_{0,1}^2}{|M|}\,,
\label{E-ple-12}\\[5pt]
\frac{|M|}{(1-\varepsilon)^2\pi j_{0,1}^2}\, \lambda_k(M) \ge k &\ge
\frac{|M|}{4 \pi}\, \lambda_k(M) - A(M)\, \sqrt{\lambda_k(M)} - B(M) +1\,,
\label{E-ple-14}
\end{align}
In particular, $F_{M,\varepsilon}\left(\lambda_k(M)\right) \le 0$, where
\begin{equation}\label{E-ple-16}
F_{M,\varepsilon}(\lambda) = \frac{|M|}{4\pi}\left( 1 - \frac{4}{(1-\varepsilon)^2j_{0,1}^2} \right)\, \lambda - A(M) \sqrt{\lambda} - B(M) +1\,,
\end{equation}
so that $\sqrt{\lambda_k(M)} \le D(M,\varepsilon)$, where $D(M,\varepsilon)$ is the largest root of the equation $F_{M,\varepsilon}(\lambda)=0$.
\end{lemma}%

\pf To prove \eqref{E-ple-12}, choose any eigenfunction associated with $\lambda_k(M)$ and having precisely $k$ nodal domains.  Since $k \ge \frac{1}{C(M,\varepsilon)}$, one of them, call it $\omega$, has area $|\omega| \le \frac{|M|}{k} \le C(M,\varepsilon)\, |M|$, so that its first Dirichlet eigenvalue satisfies  \eqref{E-ple-2}. Then, use the fact that $\lambda_k(M) = \delta_1(\omega)$.\smallskip

To prove \eqref{E-ple-14}, use the fact that $k-1 = W_M\left( \lambda_k(M)\right)$ (this is because the eigenvalue is Courant-sharp), and apply inequalities \eqref{E-ple-4} and \eqref{E-ple-12}. \smallskip

Since $j_{0,1}>2$, choosing $\varepsilon < 1 - \frac{2}{j_{0,1}}$, the coefficient of the leading term in $F_{M,\varepsilon}$ is positive, and the function tends to infinity when $\lambda$ tends to infinity.\hfill \qed

\begin{corollary}\label{C-ple-2}
For $\varepsilon < 1 - \frac{2}{j_{0,1}}$, in order to be Courant-sharp, the eigenvalue $\lambda_k(M)$ must satisfy $\lambda_{k-1}(M) < \lambda_k(M)$, and
\begin{equation}\label{E-ple-18}
\left\{
\begin{array}{l}
\text{either~} k < \frac{1} {C(M,\varepsilon)} \,,\\[5pt]
\text{or~} k \ge \frac{1}{C(M,\varepsilon)}, ~\sqrt{\lambda_k(M)} \le D(M,\varepsilon), \text{~and~~} \frac{\lambda_k(M)}{k} \ge (1-\varepsilon)^2 \frac{\pi j_{0,1}^2}{|M|}\,.
\end{array}%
\right.
\end{equation}
To conclude whether the eigenvalue $\lambda_k(M)$ is actually Courant-sharp, it remains to determine the maximum number of nodal domains of an eigenfunction in the eigenspace $\cE\left( \lambda_k(M) \right)$.
\end{corollary}%

\begin{remarks}\label{R-ple-2}~\vspace{-2mm}
\begin{enumerate}
  \item In Sections~\ref{S-kle} and \ref{S-cyl}, we will use the fact that we can choose $\varepsilon=0$ in the isoperimetric inequality \eqref{E-ple-2} for the flat Klein bottles and for the flat cylinders. This follows from \cite{HHM}, and was already used in \cite{Le,BH2016,BHK2020} for the flat torus and for the M\"{o}bius band.
  \item When $\partial M \not = \emptyset$, and for Neumann or Robin eigenfunctions of $(M,g)$, the inequality \eqref{E-ple-2} can only be applied to a nodal domain which does not touch the boundary $\partial M$. To prove a result \`{a} la Pleijel for Neumann or Robin eigenfunctions, as in \cite{Pol2009, Le2019,HP2015,GH1,GH2}, it is necessary to take care of this difficulty.
  \item The proof of Lemma~\ref{L-ple-2} also yields the following inequality.
   Let $\kappa(k)$ denote the maximal number of nodal domains of an
   eigenfunction associated with $\lambda_k$. Then,
$$
\limsup_{k} \frac{\kappa(k)}{k} \le \gamma(2) := \left( \frac{2}{j_{0,1}} \right)^2 < 1.
$$
This inequality, which generalizes Pleijel's inequality \cite{Pl} valid for plane domains, is actually a particular case of a general result valid for any closed Riemannian manifold, \cite{BM}. The interesting feature is that the upper bound $\gamma(n) < 1$ only depends on the dimension.
\end{enumerate}
\end{remarks}%

\section{Preliminaries on Klein bottles}\label{S-pre}

\subsection{Klein bottles}

In this note, we are interested in the flat Klein bottles. More precisely, given $a,b > 0$, we consider the isometries of $\R^2$ given by
\begin{equation}\label{E-pre-4}
\left\{
\begin{array}{l}
\tau_1: (x,y) \mapsto (x,y+b)\,,\\[5pt]
\tau_2: (x,y) \mapsto (x+a,y)\,,\\[5pt]
\tau: (x,y) \mapsto (x+\frac{a}{2},b-y)\,.
\end{array}%
\right.
\end{equation}

We denote by $G_2$ (resp. $G$) the group generated by $\tau_1$ and $\tau_2$ (resp. by $\tau_1$ and $\tau$). These groups act properly and freely by isometries on $\R^2$ equipped with the usual scalar product. Since $\tau^2 = \tau_2$, the group $G_2$ is a subgroup of index $2$ of the group $G$. We denote by $\T_{a,b}$ (resp. $\K_{a,b}$) the torus $\R^2/G_2$ (resp. the Klein bottle $\R^2/G$). We equip $\T_{a,b}$ and $\K_{a,b}$ with the induced flat Riemannian metrics. \medskip

A fundamental domain for the action of $G_2$ (resp. $G$) on $\R^2$ is the rectangle $\cT_{a,b} = (0,a)\times (0,b)$ (resp. the rectangle $\cK_{a,b} = (0,\frac{a}{2})\times (0,b)$, see Figure~\ref{F-pre-4}\,(A)). The horizontal sides of $\cK_{a,b}$ are identified with the same orientation, the vertical sides are identified with the opposite orientations.\medskip

The geodesics of the Klein bottle are the images of the lines in $\R^2$ under the  Riemannian covering map $\R^2 \to \K_{a,b}$ (see \cite{GHL}). They can be looked at in the fundamental domain $\cK_{a,b}$, taking into account the identifications $(x,0) \sim (x,b)$ and $(0,y) \sim (\frac a2,b-y)$. Among them, we have some special geodesics, see Figure~\ref{F-pre-4}\,(B)-(C),
\begin{itemize}
  \item[$\diamond$] $t \mapsto (t,0)$ and $t \mapsto (t,\frac b2)$,  for $0 \le t \le \frac a2$, which are periodic geodesics of length $\frac a2$;
  \item[$\diamond$] for $0 < y_0 < \frac b2$, $\gamma_{y_0}: t \mapsto \left\{
      \begin{array}{ll}
      (t,y_0), &0\le t \le \frac a2,\\[5pt]
      (t-\frac a2,b-y_0),&\frac a2 \le t \le a,
      \end{array}
      \right.$\\
      which is a periodic geodesic of length $a$; the two horizontal lines in blue in Figure~\ref{F-pre-4}\,(C) yield a periodic geodesic of the Klein bottle;
  \item[$\diamond$] for $0 \le x_0 \le \frac a2$, $t \mapsto (x_0,t)$,  with $0 \le t \le b$, is a periodic geodesic of length $b$.
\end{itemize}\medskip

\begin{remark}\label{R-pre-2}
The description of geodesics of the Klein bottles as projected lines implies that the shortest, nontrivial, periodic geodesic of $\K_{a,b}$ has length $\min\set{\frac a2,b}$.
\end{remark}%

\begin{remark}\label{R-pre-4}
Scissoring the Klein bottle along the  blue geodesic $t \mapsto \gamma_{y_0}(t)$, with $0 < y_0 < \frac b2$ and $0\le t \le a$, divides the surface into two M\"{o}bius strips whose center lines are the geodesics $t \mapsto (t,0)$ and $t \mapsto (t,\frac b2)$, see Figure~\ref{F-pre-4}\,(D).
\end{remark}%

The isometry $\tau$ of $\R^2$ induces an isometry on the torus $\T_{a,b}$ so that we can identify $\K_{a,b}$ with the quotient $\T_{a,b}/\set{Id,\tau}$. It follows that the eigenfunctions of the Klein bottle $\K_{a,b}$ are precisely the eigenfunctions of the torus $\T_{a,b}$ which are invariant under the map $\tau$. Because $\tau$ is orientation reversing, the surface $\K_{a,b}$ is non-orientable with orientation double cover $\T_{a,b}$.

\newcommand{\scaK}{0.33}
\begin{figure}[htb]
\centering
\begin{subfigure}[t]{\scaK\textwidth}
\centering
\includegraphics[width=\linewidth]{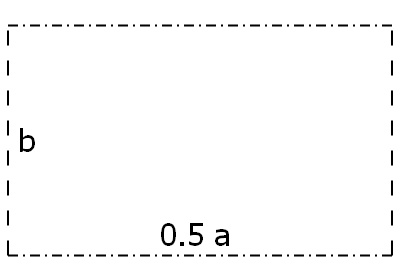}
\caption{}
\end{subfigure}
\begin{subfigure}[t]{\scaK\textwidth}
\centering
\includegraphics[width=\linewidth]{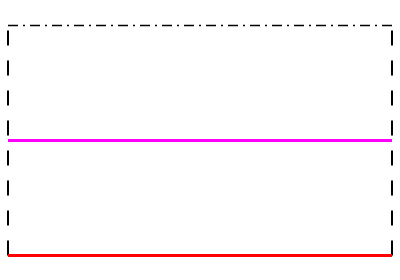}
\caption{}
\end{subfigure}
\begin{subfigure}[t]{\scaK\textwidth}
\centering
\includegraphics[width=\linewidth]{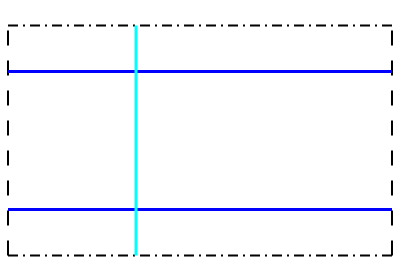}
\caption{}
\end{subfigure}
\begin{subfigure}[t]{\scaK\textwidth}
\centering
\includegraphics[width=\linewidth]{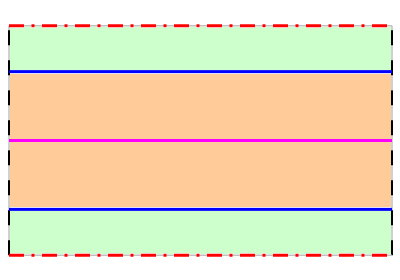}
\caption{}
\end{subfigure}
\caption{Fundamental domain, geodesics, partition into M\"{o}bius strips}\label{F-pre-4}
\end{figure}

\FloatBarrier

\subsection{The spectrum of Klein bottles}

A complete family of (complex) eigenfunctions of the flat torus $\T_{a,b}$ is
\begin{equation}\label{E-pre-10}
f_{m,n}(x,y) = \exp\left( i\frac{2\pi mx}{a} \right) \, \exp\left( i\frac{2\pi ny}{b}\right)\,, ~~m,n \in \Z\,,
\end{equation}
with associated eigenvalues $\lambdah(m,n) = 4\pi^2 \left( \frac{m^2}{a^2} + \frac{n^2}{b^2}\right)$. Given some eigenvalue $\lambda$ of $\T_{a,b}$, we introduce the set,
\begin{equation}\label{E-pre-12}
\cL_{\lambda} = \set{(m,n) \in \Z^2 ~|~ \lambdah(m,n) = \lambda}
\end{equation}

A general (complex) eigenfunction of $\T_{a,b}$, with eigenvalue $\lambda$, is of the form
\begin{equation}\label{E-pre-14}
\phi = \sum_{(m,n) \in \cL_{\lambda}} \alpha_{m,n}\, f_{m,n}\,,
\end{equation}
with $\alpha_{m,n} \in \C$. \\
The function $\phi$ in invariant under $\tau$, $\phi = \phi\circ \tau$, if and only if
\begin{equation*}
\sum_{(m,n) \in \cL_{\lambda}} \alpha_{m,n}\, f_{m,n}(x,y) = \sum_{(m,n) \in \cL_{\lambda}} \alpha_{m,n} (-1)^m \,f_{m,-n}(x,y)\,,
\end{equation*}
or, equivalently, if and only if
\begin{equation}\label{E-pre-16}
\alpha_{m,-n} = (-1)^m \alpha_{m,n}\,, ~~\forall (m,n) \in \cL_{\lambda}\,.
\end{equation}

We can rewrite a $\tau$-invariant eigenfunction $\phi$ as,
\begin{equation}\label{E-pre-18}
\phi = \sum_{(m,0) \in \cL_{\lambda}, m \text{~even}} \alpha_{m,0}\, f_{m,0} + \sum_{(m,n) \in \cL_{\lambda}, n > 0} \alpha_{m,n} \, \left( f_{m,n} + (-1)^m f_{m,-n}\right)\,.
\end{equation}

The following lemma follows readily.

\begin{lemma}[\cite{BGM}]\label{L-pre-12}
A complete family of real eigenfunctions of the flat Klein bottle $\K_{a,b}$ is given by the following functions.
\begin{equation}\label{E-pre-20}
\left\{
\begin{array}{ll}
\text{For~}  m = 0, n\in \N :& \cos\left( \frac{2\pi ny}{b}\right);\\[8pt]
\text{for~} m\in \Nb \text{~even}, n \in \N :& \cos\left( \frac{2\pi mx}{a}\right)\cos\left( \frac{2\pi ny}{b}\right);\, \sin\left( \frac{2\pi mx}{a}\right)\cos\left( \frac{2\pi ny}{b}\right);\\[8pt]
\text{for~} m \in \Nb \text{~odd}, n \in \Nb :& \cos\left( \frac{2\pi mx}{a}\right)\sin\left( \frac{2\pi ny}{b}\right);\, \sin\left( \frac{2\pi mx}{a}\right)\sin\left( \frac{2\pi ny}{b}\right).
\end{array}%
\right.
\end{equation}
Here, $\N$ denotes the set of non-negative integers, and $\Nb$ the set of positive integers.
\end{lemma}%

\begin{remark}\label{R-pre-6}
If $\cL_{\lambda} \cap ( \set{0}\times \Z) \not = \emptyset$, the multiplicity of $\lambda$ is odd; if $\cL_{\lambda} \cap (\set{0}\times \Z) = \emptyset$, the multiplicity of $\lambda$ is even.
\end{remark}%

\subsection{Choices for $a$ and $b$}

In this paper, we restrict our attention to the case $a=b=2\pi$, i.e. to the flat Klein bottle $\K_1 := \K_{2\pi,2\pi}$, whose fundamental domain is the rectangle $\cK_1 = (0,\pi) \times (0,2\pi)$, and  to the case $a=2\pi$, $b=\pi$, i.e., to the flat Klein bottle $\K_2 := \K_{2\pi,\pi}$, whose fundamental domain is the square $\cK_{2} = (0,\pi)\times(0,\pi)$. As in \cite{HHO3} for flat tori, we could consider other values of the pair $(a,b)$, in particular we could look at what happens when $a$ is fixed, $b$ tends to zero, and vice-versa.\medskip

We denote the associated square flat tori by $\T_1$ and $\T_2$ respectively, with corresponding fundamental domains $\cT_1 = (0,2\pi) \times (0,2\pi)$ and $\cT_2 = (0,2\pi)\times (0,\pi)$.\medskip

As points of the spectrum, the eigenvalues of the flat Klein bottle $\K_c, c \in \set{1,2}$ are the numbers of the form $\lambdah(p,q)=p^2+c^2q^2$,  with $p,q \in \N$, and the extra condition that $p$ is even when $q=0$. As usual, the eigenvalues of $\K_c$ are listed in nondecreasing order, multiplicities accounted for, starting from the label $1$,
\begin{equation}\label{E-pre-32}
0 = \lambda_1(\K_c) < \lambda_2(\K_c) \le \lambda_3(\K_c) \le \cdots .
\end{equation}

For $\lambda \ge 0$, we introduce the \emph{Weyl counting function},
\begin{equation}\label{E-pre-34}
W_{\K_c}(\lambda) = \# \set{j ~|~ \lambda_j(\K_c) < \lambda} \,.
\end{equation}

Weyl's asymptotic law tells us that
\begin{equation}\label{E-pre-36}
W(\lambda) = \frac{|\K_c|}{4\pi}\, \lambda + O(\sqrt{\lambda}) = \frac{\pi}{2c}\, \lambda  + O(\sqrt{\lambda}),
\end{equation}
where $|\K_c|$ denotes the area of the Klein bottle, namely $|\K_c| = \frac{2 \pi^2}{c}$.\medskip

For later purposes, we also introduce the set
\begin{equation}\label{E-pre-38}
\cL_c(\lambda) := \set{(m,n) \in \N^2 ~|~ m^2+c^2 n^2 < \lambda}
\end{equation}
and the \emph{counting function},
\begin{equation}\label{E-pre-40}
L_c(\lambda) = \#\left(\set{(m,n) \in \N^2 ~|~ m^2+c^2 n^2 < \lambda}\right)\,.
\end{equation}

\section{Courant-sharp eigenvalues of the Klein bottles $\K_c, c\in \set{1,2}$}\label{S-kle}

The purpose of this section is to determine the Courant-sharp eigenvalues of the Klein bottles $\K_c$, following Pleijel's method, Section~\ref{S-ple}.

\begin{lemma}\label{L-kle-2}
Let $\K_c, c\in \set{1,2}$ be the flat Klein bottles introduced in Section~\ref{S-pre}. Let $\lambda_k(\K_c)$ be a Courant-sharp eigenvalue of $\K_c$, with label $k \ge \frac{2\pi}{c}$ (i.e. $k \ge 7$ when $c=1$; $k\ge 4$ when $c=2$).\smallskip

When $c=1$,
\begin{align}
\frac{\lambda_k(\K_1)}{k} &\ge \frac{j_{0,1}^2}{2\pi} \ge 0.920422\,,
\label{E-kle-2}\\
\frac{2\pi}{j_{0,1}^2}\lambda_k(\K_1) \ge k &\ge \frac{\pi}{2}\lambda_k(\K_1) - 2 \sqrt{\lambda_k(\K_1)} -2\,.\label{E-kle-4}
\end{align}

When $c=2$,
\begin{align}
\frac{\lambda_k(\K_2)}{k} &\ge \frac{j_{0,1}^2}{\pi} \ge 1.840844\,,
\label{E-kle-6}\\
\frac{\pi}{j_{0,1}^2}\lambda_k(\K_2) \ge k &\ge \frac{\pi}{4}\lambda_k(\K_1) - \frac 32\sqrt{\lambda_k(\K_2)} -1\,.\label{E-kle-8}
\end{align}

In particular, $\lambda_k(\K_1) < 25$ and $\lambda_k(\K_2) \le 47$.
\end{lemma}%

\pf By Remark~\ref{R-pre-2}, the shortest closed geodesic of $\K_c, c \in \set{1,2}$, has length $\pi$. According to \cite[\S~7]{HHM}, any domain $\omega \subset \K_c$, with area $|\omega| \le \pi$, satisfies the Euclidean isoperimetric inequality $|\partial \omega|^2 \ge 4\pi |\omega|$ (where $|\partial \omega|$ denotes the length of the boundary $\partial \omega$). It follows that
\begin{equation*}
\delta_1(\omega) \ge \frac{\pi j_{0,1}^2}{|\omega|}.
\end{equation*}

If $\lambda_k(\K_c)$ is Courant-sharp with $k \ge \frac{2\pi}{c}$, then, for any associated eigenfunction $u$ with precisely $k$ nodal domains, there exists at least one nodal domain $\omega$ with $|\omega| \le \frac{\K_c}{k} \le \pi$, and we have
\begin{equation*}
\lambda_k(\K_c) = \delta_1(\omega) \ge \frac{\pi j_{0,1}^2}{|\omega|} \ge k \frac{c j_{0,1}^2}{2\pi}\,.
\end{equation*}

This proves inequalities~\eqref{E-kle-2} and \eqref{E-kle-6}.\smallskip

To prove the other inequalities, consider the set
\begin{equation*}
\cE_c(\lambda) := \set{(x,y) ~|~ 0 \le x, y, ~ x^2+c^2y^2 < \lambda}.
\end{equation*}

Then,
\begin{equation*}
\cE_c(\lambda)   \subset \bigcup_{(m,n) \in \cL_c(\lambda)}\, [m,m+1]\times [n,n+1]\,,
\end{equation*}
and hence
\begin{equation}\label{E-kle-10}
L_c(\lambda) = \#\left( \cL_c(\lambda) \right) \ge |\cE_c(\lambda)| = \frac{\pi}{4c}\, \lambda\,.
\end{equation}

From the description of the spectrum of $\K_c$ (Lemma~\ref{L-pre-12}), the contribution of a pair $(m,n) \in \cL_c(\lambda)$, to $W_{\K_c}(\lambda)$ is
 \begin{enumerate}
   \item[a)] $2$ if $m, n \ge 1$;
   \item[b)] $1$ if $m=0$, and $n \ge 0$;
   \item[c)] $2$ if $m\ge 2$ is even, and $n=0$;
   \item[d)] $0$ otherwise.
 \end{enumerate}

Recall that multiplicities also arise from the number of solutions of the equation $m^2+c^2n^2 = m_0^2+c^2n_0^2$. It follows that,
\begin{align*}
W_{\K_c}(\lambda) &= 2\, \left( L_c(\lambda) - \lfloor\sqrt{\lambda}\rfloor - \lfloor\frac{\sqrt{\lambda}}{c}\rfloor -1 \right) + \left( \lfloor\frac{\sqrt{\lambda}}{c}\rfloor +1\right) + 2\lfloor\frac{\sqrt{\lambda}}{2}\rfloor\,,\\[5pt]
&= 2\, L_c(\lambda) - 2\, \lfloor\sqrt{\lambda}\rfloor - \lfloor\frac{\sqrt{\lambda}}{c}\rfloor + 2\, \lfloor\frac{\sqrt{\lambda}}{2}\rfloor -1\,,
\end{align*}
where $\lfloor x\rfloor$ denotes the integer part of $x\ge 0$. It follows that, for all $\lambda \ge 0$,
\begin{equation*}
W_{\K_c}(\lambda) \ge \left\{
\begin{array}{l}
\frac{\pi}{2}\lambda - 2 \sqrt{\lambda} -3 \hspace{2mm} \text{when~} c=1\,,\\[5pt]
\frac{\pi}{4}\lambda - \frac 32 \sqrt{\lambda} -2 \hspace{2mm} \text{when~} c=2\,.
\end{array}
\right.
\end{equation*}

Inequalities~\eqref{E-kle-4} and \eqref{E-kle-8} follow from the previous inequalities. \hfill \qed \medskip

\begin{table}[ht!]
\resizebox{0.95\textwidth}{!}{%
  \begin{minipage}[t]{0.5\linewidth}
    \centering
  \begin{tabular}[c]{|c|c|c|c|}%
\hline &&&\\[2pt]
$\K_1$ & $\lambda_{k_{\min}}(\K_1)$ & $\lambda_{k_{\max}}(\K_1)$ &
$\dfrac{\lambda_{k_{\min}}(\K_1)}{k_{\min}}$\\[8pt]
\hline
$1$ & $\lambda_2$ & $\lambda_2$ & --- \\[5pt]
\hline
$2$ & $\lambda_3$ & $\lambda_4$ & --- \\[5pt]
\hline
$4$ & $\lambda_5$ & $\lambda_7$ & --- \\[5pt]
\hline
$5$ & $\lambda_8$ & $\lambda_{11}$ & $0.625$ \\[5pt]
\hline
$8$ & $\lambda_{12}$ & $\lambda_{13}$ & $0.6666$ \\[5pt]
\hline
$9$ & $\lambda_{14}$ & $\lambda_{14}$ & $0.6428$ \\[5pt]
\hline
$10$ & $\lambda_{15}$ & $\lambda_{18}$ & $0.6666$ \\[5pt]
\hline
$13$ & $\lambda_{19}$ & $\lambda_{22}$ & $0.6842$ \\[5pt]
\hline
$16$ & $\lambda_{23}$ & $\lambda_{25}$ & $0.6956$ \\[5pt]
\hline
$17$ & $\lambda_{26}$ & $\lambda_{29}$ & $0.6538$ \\[5pt]
\hline
$18$ & $\lambda_{30}$ & $\lambda_{31}$ & $0.6$ \\[5pt]
\hline
$20$ & $\lambda_{32}$ & $\lambda_{35}$ & $0.625$ \\[5pt]
\hline
$25$ & $\lambda_{36}$ & $\lambda_{40}$ & $0.6944$ \\[5pt]
\hline
\end{tabular}
  \end{minipage}%
  \begin{minipage}[t]{0.5\linewidth}
    \centering
\begin{tabular}[c]{|c|c|c|c|}%
\hline &&&\\[2pt]
$\K_2$ & $\lambda_{k_{\min}}(\K_2)$ & $\lambda_{k_{\max}}(\K_2)$ &
$\dfrac{\lambda_{k_{\min}}(\K_2)}{k_{\min}}$\\[8pt]
\hline
$4$ & $\lambda_2$ & $\lambda_4$ & --- \\[5pt]
\hline
$5$ & $\lambda_5$ & $\lambda_6$ & 1 \\[5pt]
\hline
$8$ & $\lambda_7$ & $\lambda_{8}$ & $1.1429$ \\[5pt]
\hline
$13$ & $\lambda_{9}$ & $\lambda_{10}$ & $1.4444$ \\[5pt]
\hline
$16$ & $\lambda_{11}$ & $\lambda_{13}$ & $1.4545$ \\[5pt]
\hline
$17$ & $\lambda_{14}$ & $\lambda_{15}$ & $1.2143$ \\[5pt]
\hline
$20$ & $\lambda_{16}$ & $\lambda_{19}$ & $1.2500$ \\[5pt]
\hline
$25$ & $\lambda_{20}$ & $\lambda_{21}$ & $1.2500$ \\[5pt]
\hline
$29$ & $\lambda_{22}$ & $\lambda_{23}$ & $1.3182$ \\[5pt]
\hline
$32$ & $\lambda_{24}$ & $\lambda_{25}$ & $1.3333$ \\[5pt]
\hline
$36$ & $\lambda_{26}$ & $\lambda_{28}$ & $1.3846$ \\[5pt]
\hline
$37$ & $\lambda_{29}$ & $\lambda_{30}$ & $1.2759$ \\[5pt]
\hline
$40$ & $\lambda_{31}$ & $\lambda_{34}$ & $1.2903$ \\[5pt]
\hline
$41$ & $\lambda_{35}$ & $\lambda_{36}$ & $1.1714$ \\[5pt]
\hline
$45$ & $\lambda_{37}$ & $\lambda_{38}$ & $1.2162$ \\[5pt]
\hline
$52$ & $\lambda_{39}$ & $\lambda_{42}$ & $1.3333$ \\[5pt]
\hline
\end{tabular}
\end{minipage}
}
\par\vspace{8pt} 
\caption{The first eigenvalues of $\K_1$ (left) and $\K_2$ (right)}\label{Tab-kle-2}
\end{table}

Table~\ref{Tab-kle-2} displays the eigenvalues of $\K_c$ less than or equal to $25$ ($c=1$), resp. $47$ ($c=2$), the corresponding labeled eigenvalues, and the ratio $\dfrac{\lambda_{k_{\min}}(\K_c)}{k_{\min}}$ which should be larger than or equal to
\begin{align*}
  \frac{j_{0,1}^2}{2\pi} \approx 0.92042 & \hspace{5mm }\text{~if~} \lambda_{k_{\min}}(\K_1) \text{~is Courant-sharp}, \\[5pt]
  \frac{j_{0,1}^2}{\pi} \approx 1.84084 & \hspace{5mm }\text{~if~} \lambda_{k_{\min}}(\K_2) \text{~is Courant-sharp}.
\end{align*}

Since the eigenvalues $\lambda_1(\K_c), \lambda_2(\K_c)$ are Courant-sharp, we conclude from Table~\ref{Tab-kle-2} that Theorem~\ref{T-int-2K} is proved in the case $c=2$. To finish the proof of the theorem when $c=1$, it remains to investigate  $\lambda_3(\K_1)$ and $\lambda_5(\K_1)$. This is done in the next section. \vspace{1pt}

\section{The eigenvalues $\lambda_3(\K_1)$ and $\lambda_5(\K_1)$ are not Courant-sharp}\label{S-smal}

In order to finish the proof of Theorem~\ref{T-int-2K} for the Klein bottle $\K_1$, we investigate the eigenvalues $\lambda_3(\K_1) = \lambdah(1,1)$ and $\lambda_5(\K_2)=\lambdah(2,0)=\lambdah(0,2)$.

\subsection{The eigenvalue $\lambda_3(\K_1)$ is not Courant-sharp}

\begin{figure}[ht]
  \centering
  \includegraphics[scale=0.3]{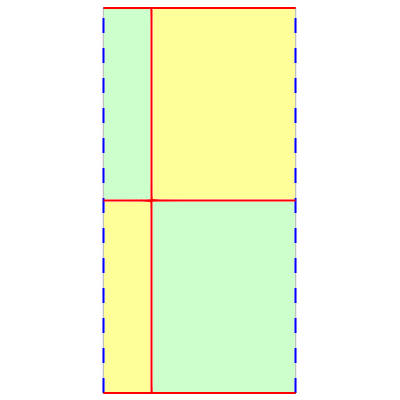}\\
  \caption{Nodal domains of $\sin(x-\alpha) \sin(y)$}\label{F-lam3}
\end{figure}

A general eigenfunction associated with $\lambda_3$ has the form $\left( A \cos(x) + B\sin(x)\right) \, \sin(y)$. It is sufficient to look at eigenfunctions of the form $\sin(x-\alpha)\, \sin(y)$. These eigenfunctions have exactly two nodal domains in $\K_1$. It follows that $\lambda_3$ is not Courant-sharp, see Figure~\ref{F-lam3}.\medskip

\subsection{The eigenvalue $\lambda_5(\K_1)$ is not Courant-sharp}

A general eigenfunction associated with $\lambda_5$ has the form $A \cos(2x) + B \sin(2x) + C \cos(2y)$. Up to multiplication by a scalar, it suffices to consider the family $\cos\theta \cos(2x-\alpha) + \sin\theta \cos(2y)$, with $\theta \in [0,\pi)$ and $\alpha \in [0,\pi)$. Choosing the fundamental domain appropriately, we can assume that $\alpha = 0$. Changing $y$ to $y+\pib$, we see that it suffices to consider $\theta \in [0,\pib]$. The nodal sets are known explicitly when $\theta = 0$ or $\pib$. We now consider the family
\begin{equation}\label{E-smal-4}
\phi_{\theta}(x,y) = \cos\theta \cos(2x) + \sin\theta \cos(2y), ~~\theta \in (0,\pib).
\end{equation}

The critical zeros (points at which both the function and its differential vanish) satisfy the system,
\begin{equation}\label{E-small-6}
\left\{
\begin{array}{l}
\cos\theta \cos(2x) + \sin\theta \cos(2y) = 0\,,\\[5pt]
\sin(2x)=0\,,\\[5pt]
\sin(2y) = 0\,.
\end{array}%
\right.
\end{equation}

\begin{figure}[h]
\centering
\begin{subfigure}[t]{.30\textwidth}
\centering
\includegraphics[width=\linewidth]{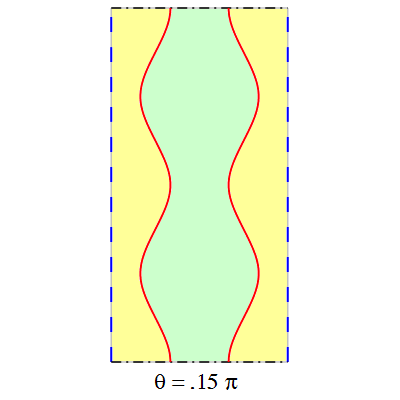}
\caption{}
\end{subfigure}
\begin{subfigure}[t]{.30\textwidth}
\centering
\includegraphics[width=\linewidth]{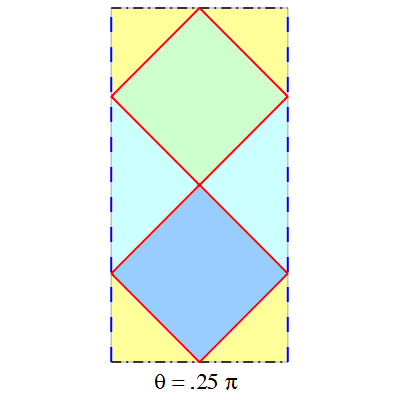}
\caption{}
\end{subfigure}
\begin{subfigure}[t]{.30\textwidth}
\centering
\includegraphics[width=\linewidth]{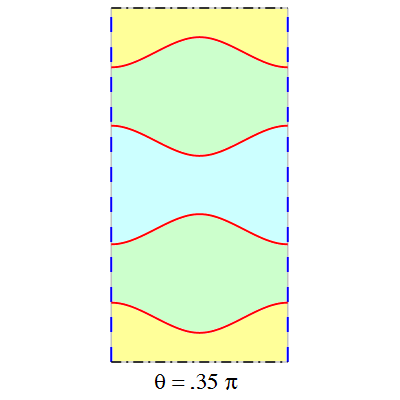}
\caption{}
\end{subfigure}
\caption{Nodal sets for $\cos\theta \cos(2x) + \sin\theta \cos(2y)$}\label{F-ev5-c}
\end{figure}

It follows that critical zeros only occur for $\theta = \pid$, so that the nodal set of $\phi_{\theta}$ is a regular curve when $\theta \not = \pid$. The nodal set of the eigenfunction
\begin{equation}\label{E-smal-8}
\phi_{\pid}(x,y) = \cos(2x) + \cos(2y) = 2 \cos(x+y) \cos(x-y)
\end{equation}
is explicit, see Figure~\ref{F-ev5-c}\,(B). An analysis \`{a} la Stern, see \cite[Section~5]{BHK2020} or \cite{BH-s},  shows that the nodal sets of $\phi_{\theta}$ for $0 < \theta < \pid$ and $\pid < \theta < \pib$ are given by Figures~\ref{F-ev5-c}\,(A) and (C) respectively. More precisely, we first note that the common zeros to $\phi_{0}$ and $\phi_{\pib}$ are common zeros to all $\phi_{\theta}$. Since $\theta \in (0,\pib)$, and hence $\sin\theta \, \cos\theta > 0$, it follows that, except for the common zeros, the nodal set of $\phi_{\theta}$ is contained in the set $\cos(2x) \, \cos(2y) < 0$. Finally, depending on the sign of $\pid - \theta$, we can use the nodal set of $\cos(2x)$ or the nodal set of $\cos(2y)$ as barrier to obtain the behaviour described in the figures.\medskip

It follows that an eigenfunction associated with $\lambda_5$ has at most $4$ nodal domains in $\K_1$, and hence that $\lambda_5(\K_1)$ is not Courant-sharp.\medskip

This completes the proof of Theorem~\ref{T-int-2K} in the case $c=1$.

\section{Courant-sharp Dirichlet eigenvalues of flat cylinders}\label{S-cyl}

\subsection{Preliminaries}\label{SS-cyl-1}

Given $r > 0$, let $\cC_r = (0,\pi) \times \bS^1_r$ denote the cylinder with radius $r$ and length $\pi$, equipped with the flat product metric. We consider the Dirichlet boundary condition on $\partial \cC_r = \set{0,\pi} \times \bS^1_r$. A complete family of Dirichlet eigenfunctions is given by,
\begin{equation}\label{E-cyl-2}
\left\{
\begin{array}{l}
\sin(mx) \, \cos\left( \frac{ny}{r}\right)\,, ~\text{~for~} m\in \Nb, n\in \N\,,\\[5pt]
\sin(mx) \, \sin\left( \frac{ny}{r}\right)\,, ~\text{~for~} m\in \Nb, n\in \Nb\,,
\end{array}%
\right.
\end{equation}
with associated eigenvalues the numbers $\lambdah_r(m,n) = m^2 + \frac{n^2}{r^2}$.\medskip

For $m \in \Nb$, the point $(m,0)$ contributes one eigenfunction $\sin(mx)$; for $m, n \in \Nb$, the point $(m,n)$ contributes two functions, $\sin(mx) \, \cos\left( \frac{ny}{r}\right)$ and $\sin(mx) \, \sin\left( \frac{ny}{r}\right)$. Multiplicities of the eigenvalues also occur when the equation $m^2 + \frac{n^2}{r^2}= m_0^2 + \frac{n_0^2}{r^2}$ has multiple solutions within the above range.\medskip

As usual, we arrange the Dirichlet eigenvalues of $\cC_r$ in non-decreasing order, starting from the label $1$, multiplicities accounted for,
\begin{equation}\label{E-cyl-4}
0 < \lambda_1(\cC_r) < \lambda_2(\cC_r) \le \lambda_3(\cC_r) \le \cdots\,.
\end{equation}

The purpose of this section is to prove Theorem~\ref{T-int-2C}, i.e. to determine the Courant-sharp eigenvalues in two specific cases $r \in \set{\frac 12, 1}$ whose eigenvalues have higher multiplicities. The corresponding cylinder are the orientation covers of the M\"{o}bius band we studied in \cite{BHK2020}. As in \cite{HHO2,HHO3}, we could also consider other values of $r$, in particular $r$ close to zero or $r$ very large.\medskip

\begin{remark}\label{R-cyl-2}
As pointed out in the introduction, $\lambda_1(\cC_r)$ and $\lambda_2(\cC_r)$ are always Courant-sharp, and any Courant-sharp eigenvalue $\lambda_k(\cC_r)$ satisfies the inequality $\lambda_k(\cC_r) > \lambda_{k-1}(\cC_r)$.
\end{remark}%

Following Pleijel's method, Section~\ref{S-ple}, we introduce Weyl's counting function,
\begin{equation}\label{E-cyl-10}
W_{\cC_r}(\lambda) = \# \set{j ~|~ \lambda_j(\cC_r) < \lambda}\,.
\end{equation}

According to Weyl's law,
\begin{equation}\label{E-cyl-12}
W_{\cC_r}(\lambda) = \frac{|\cC_r|}{4\pi}\, \lambda + O(\sqrt{\lambda}) = \frac{r\pi}{2}\, \lambda + O(\sqrt{\lambda})\,,
\end{equation}
where $|\cC_r|$ denotes the area of the cylinder. \medskip

\subsection{Courant-sharp eigenvalues of $\cC_r, r \in \set{\frac 12, 1}$}\label{SS-cyl-2}

According to \cite[\S~6]{HHM}, a domain $\omega \subset \cC_r$, with area $|\omega| \le 4\pi r^2$ satisfies a Euclidean isoperimetric inequality, and hence its least Dirichlet eigenvalue satisfies
\begin{equation}\label{E-cyl-20}
\delta_1(\omega) \ge \frac{\pi j_{0,1}^2}{|\omega|}.
\end{equation}

\begin{lemma}\label{L-cyl-2}
Let $r \in \set{\frac 12,1}$. Let $\lambda_k(\cC_r)$ be a Courant-sharp eigenvalue of $\cC_r$, with $k \ge \frac{\pi}{2r}$. The following inequalities hold.
\begin{equation}\label{E-cyl-22}
\frac{\lambda_k(\cC_r)}{k} \ge \frac{j_{0,1}^2}{2\pi r}\ge \left\{
\begin{array}{l}
\frac{j_{0,1}^2}{\pi} \approx 1.840844 \hspace{5mm }\text{~if~} r = \frac 12\,,\\[5pt]
\frac{j_{0,1}^2}{2\pi} \approx 0.920422 \hspace{5mm }\text{~if~} r = 1\,,
\end{array}
\right.
\end{equation}

\begin{equation}\label{E-cyl-24}
k-1 = W_{\cC_r}\left( \lambda_k(\cC_r) \right) \ge \frac{\pi r}{2} \,
\lambda_k(\cC_r) - (2r+1)\, \sqrt{\lambda_k(\cC_r)} -2\,.
\end{equation}

\begin{align}
\text{For~} r =\frac 12, & \hspace{5mm}\frac{\pi}{4}\left( 1 - \frac{4}{j_{0,1}^2}\right)\, \lambda_k(\cC_{\frac 12}) - 2 \sqrt{\lambda_k(\cC_{\frac 12})} -1 \le 0\,,
\label{E-cyl-26-1} \\[5pt]
& \hspace{5mm}\lambda_k(\cC_{\frac 12}) \le 76.25\,,\label{E-cyl-26-2}\\[5pt]
\text{For~} r =1, & \hspace{5mm}\frac{\pi}{4}\left( 1 - \frac{4}{j_{0,1}^2}\right)\, \lambda_k(\cC_1) - \frac 32 \sqrt{\lambda_k(\cC_1)} - \frac 12 \le 0\,,\label{E-cyl-26-3}\\[5pt]
& \hspace{5mm}\lambda_k(\cC_{1}) \le 42.40\,\label{E-cyl-26-4}
\end{align}
\end{lemma}%

\pf The arguments are similar to those used in the proof of Lemma~\ref{L-kle-2}. Inequality~\eqref{E-cyl-22} follows from \eqref{E-cyl-20}. The inequality in \eqref{E-cyl-24} follows from the description of the Dirichlet spectrum of $\cC_r$. The other inequalities follow from \eqref{E-cyl-24} and \eqref{E-cyl-22}. \hfill \qed \medskip

The following tables give the eigenvalues, the corresponding range of labels, and the ratios $\frac{\lambda_{k_{min}}}{k_{min}}$ for $\cC_{\frac 12}$ (left) and $\cC_1$ (right). For a Courant-sharp eigenvalue, this ratio should be greater than $1.840844$ for $\cC_{\frac 12}$, and greater than $0.920422$ for $\cC_1$.

\begin{table}[ht!]
\resizebox{0.9\textwidth}{!}{%
  \begin{minipage}[t]{0.45\linewidth}
    \centering
 \begin{tabular}{|c|c|c|c|}
    \hline
 $\lambda(\cC_{1/2})$ & $k_{min}$ & $k_{max}$ & $\frac{\lambda_{k_{min}}(\cC_{1/2})}{k_{min}}$ \\ \hline
        1  & 1 & 1 & -- \\ \hline
        4  & 2 & 2 & -- \\ \hline
        5  & 3 & 4 & -- \\ \hline
        8  & 5 & 6 & 1,60000 \\ \hline
        9  & 7 & 7 & 1,28571 \\ \hline
        13  & 8 & 9 & 1,62500 \\ \hline
        16  & 10 & 10 & 1,60000 \\ \hline
        17  & 11 & 12 & 1,54545 \\ \hline
        20  & 13 & 16 & 1,53846 \\ \hline
        25  & 17 & 19 & 1,47059 \\ \hline
        29  & 20 & 21 & 1,45000 \\ \hline
        32  & 22 & 23 & 1,45455 \\ \hline
        36  & 24 & 24 & 1,50000 \\ \hline
        37  & 25 & 26 & 1,48000 \\ \hline
        40  & 27 & 30 & 1,48148 \\ \hline
        41  & 31 & 32 & 1,32258 \\ \hline
        45  & 33 & 34 & 1,36364 \\ \hline
        49  & 35 & 35 & 1,40000 \\ \hline
        52  & 36 & 39 & 1,44444 \\ \hline
        53  & 40 & 41 & 1,32500 \\ \hline
        61  & 42 & 43 & 1,45238 \\ \hline
        64  & 44 & 44 & 1,45455 \\ \hline
        65  & 45 & 48 & 1,44444 \\ \hline
        68  & 49 & 52 & 1,38776 \\ \hline
        72  & 53 & 54 & 1,35849 \\ \hline
        73  & 55 & 56 & 1,32727 \\ \hline
        80  & 57 & 60 & 1,40351 \\ \hline
    \end{tabular}
  \end{minipage}%
  \begin{minipage}[t]{0.45\linewidth}
    \centering
 \begin{tabular}{|c|c|c|c|}
    \hline
 $\lambda(\cC_1)$ & $k_{min}$ & $k_{max}$ & $\frac{\lambda_{k_{min}}(\cC_1)}{k_{min}}$  \\ \hline
        1  & 1 & 1 & -- \\ \hline
        2  & 2 & 3 & 1,00000 \\ \hline
        4  & 4 & 4 & 1,00000 \\ \hline
        5  & 5 & 8 & 1,00000 \\ \hline
        8  & 9 & 10 & 0,88889 \\ \hline
        9  & 11 & 11 & 0,81818 \\ \hline
        10  & 12 & 15 & 0,83333 \\ \hline
        13  & 16 & 19 & 0,81250 \\ \hline
        16  & 20 & 20 & 0,80000 \\ \hline
        17  & 21 & 24 & 0,80952 \\ \hline
        18  & 25 & 26 & 0,72000 \\ \hline
        20  & 27 & 30 & 0,74074 \\ \hline
        25  & 31 & 35 & 0,80645 \\ \hline
        26  & 36 & 39 & 0,72222 \\ \hline
        29  & 40 & 43 & 0,72500 \\ \hline
        32  & 44 & 45 & 0,72727 \\ \hline
        34  & 46 & 49 & 0,73913 \\ \hline
        36  & 50 & 50 & 0,72000 \\ \hline
        37  & 51 & 54 & 0,72549 \\ \hline
        40  & 55 & 58 & 0,72727 \\ \hline
        41  & 59 & 62 & 0,69492 \\ \hline
        45  & 63 & 66 & 0,71429 \\ \hline
        49  & 67 & 67 & 0,73134 \\ \hline
        50  & 68 & 73 & 0,73529 \\ \hline
    \end{tabular}
\end{minipage}
}
\par\vspace{8pt} 
\caption{The first eigenvalues of $\cC_{\frac 12}$ (left) and $\cC_1$ (right)}\label{Tab-cyl-2}
\end{table}

According to Remark~\ref{R-cyl-2}, the eigenvalues $\lambda_1$ and $\lambda_2$ are always Courant-sharp. \medskip

The  table for $\cC_{\frac 12}$ (left) shows that it only remains to analyze the eigenvalue $5 = \lambda_3(\cC_{\frac 12}) = \lambda_4(\cC_{\frac 12})$. The associated eigenspace is generated by the eigenfunctions $\sin(x) \cos(2y)$ and $\sin(x) \sin(2y)$. Functions in this eigenspace have $2$ nodal domains. This finishes the proof of Theorem~\ref{T-int-2C} in the case $r=\frac 12$.\medskip

The table for $\cC_1$ (right) shows that it only remains to analyze two eigenvalues $4$ and $5$. The eigenvalue $4 = \lambda_4(\cC_1)$ is simple, with associated eigenfunction $\sin(2x)$ which has two nodal domains. The eigenvalue $5 = \lambda_5(\cC_1) = \cdots = \lambda_8(\cC_1)$ has multiplicity $4$. The corresponding eigenspace is generated by the eigenfunctions $\sin(2x) \cos(x), \sin(2x) \sin(x)$ and $\sin(x) \cos(2x), \sin(x) \sin(2y)$ which, according to \cite{BHK2020}, turn out to span the second eigenspace $\cE(\lambda_2(M_1))$ of the square M\"{o}bius strip, whose orientation cover is $\cC_1$. Since the eigenfunctions in $\cE(\lambda_2(M_1))$ have two nodal domains, it follows that the eigenfunctions in the eigenspace $\cE\left( \lambda_5(\cC_1)\right)$ have at most four nodal domains. This finishes the proof of Theorem~\ref{T-int-2C} in the case $r=1$. Note that one can give a precise analysis of nodal patterns in the eigenspace $\cE\left( \lambda_5(\cC_1) \right)$ by using the same arguments as in \cite[Section~4]{BHK2020}.


\vspace{1cm}

\end{document}